\documentclass[11pt]{amsart}
\usepackage{txfonts}
\usepackage{mathrsfs}
\usepackage{amsmath}
\usepackage[pagewise]{lineno}
\allowdisplaybreaks[4]
\usepackage{}
\usepackage{xypic}
\usepackage{amsfonts}
\usepackage{amssymb}
\usepackage{bbm}
\usepackage{color,xcolor}

\usepackage[colorlinks,
            linkcolor=blue,
           anchorcolor=blue,
            citecolor=blue
            ]{hyperref}



\newtheorem{thm}{Theorem}[section]
\newtheorem{cor}[thm]{Corollary}

\newtheorem{prop}[thm]{Proposition}
\newtheorem{defi}[thm]{Definition}
\newtheorem{rem}[thm]{Remark}
\numberwithin{equation}{section}

\newcommand{\be}{\begin{equation}}
\newcommand{\ee}{\end{equation}}
\newcommand{\bes}{\begin{eqnarray}}
\newcommand{\ees}{\end{eqnarray}}
\newcommand{\bess}{\begin{eqnarray*}}
\newcommand{\eess}{\end{eqnarray*}}

\newcommand{\bali}{\begin{align}}
\newcommand{\eali}{\end{align}}

\def\kk{\mathbbm{k}}

\begin{document}
\title[Higher Frobenius-Schur indicators]{Higher Frobenius-Schur indicators for semisimple Hopf algebras in positive characteristic}
\author{Zhihua Wang}
\address{Z. Wang\newline Department of Mathematics, Taizhou University,
Taizhou 225300, China}
\email{mailzhihua@126.com}
\author{Gongxiang Liu}
\address{G. Liu\newline Department of Mathematics, Nanjing University,
Nanjing 210093, China}
\email{gxliu@nju.edu.cn}
\author{Libin Li}
\address{L. Li\newline School of Mathematical Science, Yangzhou University, Yangzhou 225002, China}
\email{lbli@yzu.edu.cn}
\date{}
\subjclass[2010]{16T05}
\keywords{semisimple Hopf algebra, $S^2$-formula, Frobenius-Schur indicator, gauge invariant}

\begin{abstract} Let $H$ be a semisimple Hopf algebra over an algebraically closed field $\mathbbm{k}$ of characteristic $p>\dim_{\mathbbm{k}}(H)^{1/2}$. We show that the antipode $S$ of $H$ satisfies the equality $S^2(h)=\mathbf{u}h\mathbf{u}^{-1}$, where $h\in H$, $\mathbf{u}=S(\Lambda_{(2)})\Lambda_{(1)}$ and $\Lambda$ is a nonzero integral of $H$. The formula of $S^2$ enables us to define higher Frobenius-Schur indicators for the Hopf algebra $H$. This generalizes the notions of higher Frobenius-Schur indicators from the case of characteristic 0 to the case of characteristic $p>\dim_{\mathbbm{k}}(H)^{1/2}$. These indicators defined here share some properties with the ones defined over a field of characteristic 0. Especially, all these indicators are gauge invariants for the tensor category Rep$(H)$ of finite dimensional representations of $H$.
\end{abstract}
\maketitle
\section{\bf Introduction}
Linchenko-Montgomery \cite{LM} generalized the classical Frobenius-Schur (FS) indicators from group-theoretic result to the setting of a semisimple involutory Hopf algebra $H$.
They also defined higher FS indicators $\nu_n(V)$ by using idempotent integral $\Lambda$ of $H$, namely, \begin{equation}\label{eq}\nu_n(V)=\chi_V(\Lambda_{(1)}\cdots\Lambda_{(n)})\ \text{for}\ n\geq1,\end{equation} where $\chi_V$ is the character afforded by finite dimensional representation $V$ of $H$. The higher FS indicators were later extensively studied by Kashina-Sommerh$\ddot{a}$user-Zhu for semisimple Hopf algebras over an algebraically closed field of characteristic zero \cite{KSZ}, and by Ng-Schauenburg for semisimple quasi-Hopf algebras over the field of complex numbers \cite{NS1}. The notions of higher FS indicators
have been generalized to objects of a pivotal category \cite{NS2, NS3}.

However, the notions of  higher FS indicators for semisimple Hopf algebras over a field of positive characteristic seem not to be considered (except for those semisimple involutory Hopf algebras).
In this paper, we consider higher FS indicators for a finite dimensional semisimple Hopf algebra $H$ over an algebraically closed field $\kk$ of characteristic $p>\dim_{\mathbbm{k}}(H)^{1/2}$. We need to point out that the Hopf algebra $H$ here is not necessarily involutory unless the characteristic $p$ is larger than a certain number (see \cite{So,EG}).

For the  antipode $S$ of $H$, we first obtain a formula for $S^2$ as follows:
$$S^2(h)=\mathbf{u}h\mathbf{u}^{-1},$$
where $h\in H$, $\mathbf{u}=S(\Lambda_{(2)})\Lambda_{(1)}$ and $\Lambda$ is a nonzero integral of $H$.
According to the formula of $S^2$, we have an isomorphism of $H$-modules $$j_{\mathbf{u},V}:V\rightarrow V^{**},\ \ j_{\mathbf{u},V}(v)(f)=f(\mathbf{u}\cdot v)\ \text{for}\ v\in V, f\in V^*,$$
which is functorial in $V$. As the element $\mathbf{u}=S(\Lambda_{(2)})\Lambda_{(1)}$ is not necessarily a group-like element, the functorial isomorphism $j_{\mathbf{u}}:id\rightarrow(-)^{**}$ is not necessarily a tensor isomorphism. In other words, the category Rep$(H)$ of finite dimensional representations of $H$ is not necessarily  pivotal with respect to the structure $j_{\mathbf{u}}$. Even though, using the functorial isomorphism $j_{\mathbf{u}}$ we may still define the $n$-th FS indicator $\nu_n(V)$ of $V$ to be the trace of a certain $\kk$-linear operator as Ng-Schauenburg did in \cite{NS2}. It is similar to the case of characteristic 0 that
the $n$-th FS indicator $\nu_n(V)$ defined here can also be entirely described in terms of the integral $\Lambda$ of $H$ and the character $\chi_V$ of $H$-module $V$:
\begin{equation}\label{eqq}\nu_n(V)=\chi_V(\mathbf{u}^{-1}\Lambda_{(1)}\cdots\Lambda_{(n)})\ \text{for}\ n\geq1.\end{equation}
Moreover, the formula (\ref{eqq}) does not depend on the choice of
the nonzero integral $\Lambda$ and it recovers the original formula (\ref{eq}) when the characteristic of $\kk$ is zero and $\Lambda$ is idempotent.

Note that the formula (\ref{eqq}) can be written as $\nu_n(V)=\chi_V(\mathbf{u}^{-1}P_n(\Lambda))$ for $n\geq1$, where $P_{n}$ is the $n$-th Sweedler power map of $H$. Clearly, the $n$-th Sweedler power map $P_n$ is valid for all $n\in\mathbb{Z}$, this motivates us to extend the $n$-th FS indicator from $n\geq1$ to $n\in\mathbb{Z}$. That is, by definition, $\nu_n(V)=\chi_V(\mathbf{u}^{-1}P_n(\Lambda))$ for all $n\in\mathbb{Z}$. We find that the higher FS indicators defined over a field of characteristic $p>\dim_{\kk}(H)^{1/2}$ and the ones defined over a field of characteristic 0 share some common properties.
For instance, it is similar to the case of characteristic 0 (see \cite{KMN, KSZ}) that by replacing $V$ with the regular representation $H$, we reconstruct the $n$-th indicator of $H$, a notion defined by the trace of the map $S\circ P_{n-1}$. Also, it is similar to characteristic 0 case that $V$ and its dual $V^*$ have the same higher FS indicators. Especially, similar to the case of characteristic 0 that the $n$-th FS indicator $\nu_n(V)$ defined here is an invariant of the tensor category Rep$(H)$ for any $n\in\mathbb{Z}$ and any finite dimensional representation $V$ of $H$.

The paper is organized as follows: In Section 2, we present some basic results on semisimple Hopf algebras. In Section 3, we deduce the formula of $S^2$ by comparing two different forms of the character $\chi_H$ of the regular representation $H$.
We investigate some properties of the element $\mathbf{u}=S(\Lambda_{(2)})\Lambda_{(1)}$  and show that the integral $\Lambda$ of $H$ is cocommutative if and only if $S^2=id$.
In Section 4, we generalize the notions of higher FS indicators from characteristic 0 case to characteristic $p>\dim_{\kk}(H)^{1/2}$ case and find that the indicators defined here share some  common properties with the ones defined over a field of characteristic 0. In Section 5, we show that the $n$-th FS indicator $\nu_n(V)$ is a  gauge invariant for any integer $n$ and any finite dimensional representation $V$ of $H$.

\section{\bf Preliminaries}
Throughout this paper, $H$ is a finite dimensional semisimple Hopf algebra over an algebraically closed field $\kk$ of characteristic $p>\dim_{\mathbbm{k}}(H)^{1/2}$. We need to stress that all results presented here are also valid for the case of characteristic 0, although we only deal with the case of characteristic $p>\dim_{\mathbbm{k}}(H)^{1/2}$.

As a Hopf algebra, $H$ has a counit $\varepsilon$, antipode $S$, multiplication $m$ and comultiplication $\Delta$. The comultiplication $\Delta(a)$ will be written as $\Delta(a)=a_{(1)}\otimes a_{(2)}$ for $a\in H$, where we omit the summation sign. We denote by $\Lambda$ and $\lambda$ the left and right integrals of $H$ and $H^*$ respectively so that $\lambda(\Lambda)=1$. Since the semisimple Hopf algebra $H$ is unimodular, the left and right integrals of $H$ are the same.
We refer to \cite{Mon} for basic theory of Hopf algebras.

If $V$ is a finite dimensional $H$-module, then $V$ is also called a representation of $H$ via the algebra homomorphism $\rho_V:H\rightarrow \text{End}_{\kk}(V)$ given by $\rho_V(h)(v)=h\cdot v$ for $h\in H$ and $v\in V$. We will make no distinction between the two notions. The character of $V$ is the map $\chi_V:H\rightarrow \kk$ given by $\chi_V(h)=\text{tr}(\rho_V(h))$ for $h\in H$.
The $\kk$-linear dual space $V^*$ is also an $H$-module via $(h\cdot f)(v):=f(S(h)\cdot v)$ for $h\in H$, $f\in V^*$ and $v\in V$.
In particular, the dual module $V^*$ has the character $\chi_{V^*}=\chi_V\circ S$. The category Rep$(H)$ of finite dimensional representations of $H$ is a semisimple tensor category, where the monoidal structure stems from the comultiplication $\Delta$.

Recall that the dual Hopf algebra $H^*$ has an $H$-bimodule structure given by $$(a\rightharpoonup f)(b)=f(ba),\ (f\leftharpoonup a)(b)=f(ab)\ \text{for}\ a,b\in H,\ f\in H^*.$$ Moreover, $(H^*,\leftharpoonup)$ and $(\rightharpoonup,H^*)$ are free $H$-modules generated by $\lambda$, i.e., $H^*=\lambda\leftharpoonup H$ and $H^*=H\rightharpoonup \lambda$ (see \cite[Corollary 2(b)]{Rad}). This provides an associative and non-degenerate bilinear form $H\times H\rightarrow \kk$ by $a\times b\mapsto\lambda(ab)$ for $a,b\in H$. Moreover, the pair $(H,\lambda)$ is a Frobenius algebra with the Frobenius homomorphism $\lambda$ satisfying the equality (see \cite[Eq.(1)]{Rad}): \begin{equation}\label{equ01}a=\lambda(a\Lambda_{(1)})S(\Lambda_{(2)})=\lambda(S(\Lambda_{(2)})a)\Lambda_{(1)}\ \text{for}\ a\in H.\end{equation} The pair $\Lambda_{(1)}\otimes S(\Lambda_{(2)})$ satisfying (\ref{equ01}) is called the dual basis of $H$ with respect to the Frobenius homomorphism $\lambda$.

Since the right integral $\lambda$ of $H^*$ satisfies  $\lambda(ab)=\lambda(S^2(b)a)$ for all $a,b\in H$ (see \cite[Theorem 3(a)]{Rad}), the Hopf algebra $H$ is a symmetric algebra with a symmetric bilinear form given by  $$H\times H\rightarrow \kk,\ a\times b\mapsto\lambda(uab)=(\lambda\leftharpoonup u)(ab)=(u\rightharpoonup\lambda)(ab),$$ where $u$ is a unit of $H$ satisfying $S^2(h)=uhu^{-1}$ for all $h\in H$ and the Frobenius homomorphism $\lambda\leftharpoonup u=u\rightharpoonup\lambda$ holds because $\lambda(au)=\lambda(S^2(u)a)=\lambda(ua)$ for all $a\in H$. Using (\ref{equ01}) we may see that the pair $\Lambda_{(1)}\otimes u^{-1}S(\Lambda_{(2)})$ is a dual basis of $H$ with respect to $\lambda\leftharpoonup u\ (=u\rightharpoonup\lambda)$ (see also \cite[Lemma 1.4(2)]{Doi}). The symmetry of the Frobenius homomorphism $\lambda\leftharpoonup u\ (=u\rightharpoonup\lambda)$ means that \begin{equation}\label{equ5}\Lambda_{(1)}\otimes u^{-1}S(\Lambda_{(2)})= u^{-1}S(\Lambda_{(2)})\otimes \Lambda_{(1)}.\end{equation}

By Wedderburn's theorem, the semisimple Hopf algebra $H$ is isomorphic to a direct sum of full matrix algebras over $\kk$, namely,
$$H\cong\bigoplus_{i\in I}M_{d_i}(\kk).$$ Let $e_i$ be the idempotent of $H$ satisfying that $He_i\cong M_{d_i}(\kk)$. Then $\{e_i\}_{i\in I}$ forms a complete set of central primitive idempotents of $H$.
Let $V_i$ be a simple left module (unique up to isomorphism) over the matrix algebra $M_{d_i}(\kk)$. Then $\dim_{\kk}(V_i)=d_i$ and $\{V_i\}_{i\in I}$ forms a complete set of simple left $H$-modules up to isomorphism. The left regular representation $H$ has the decomposition $H\cong\bigoplus_{i\in I}V_i^{\oplus d_i}$ as $H$-modules, so the character $\chi_{H}$ of the left regular representation $H$ is equal to $\sum_{i\in I}d_i\chi_i$, where each $\chi_i$ is the character of $V_i$.

For any simple $H$-module $V_i$ and any $\varphi\in\text{End}_{\kk}(V_i)$, we use the dual basis $\Lambda_{(1)}\otimes u^{-1}S(\Lambda_{(2)})$ with respect to the Frobenius homomorphism $\lambda\leftharpoonup u$ to define the map $\mathcal{I}(\varphi)\in\text{End}_{\kk}(V_i)$ by $$\mathcal{I}(\varphi)(v)=\Lambda_{(1)}\varphi(u^{-1}S(\Lambda_{(2)})v)\ \text{for}\ v\in V_i.$$
Note that $\mathcal{I}(\varphi)$ lies in $\text{End}_{H}(V_i)\cong \kk$. There exists a unique element $c_{i}\in \kk$ such that
\begin{equation}\label{equ2}\mathcal{I}(\varphi)=c_i\text{tr}(\varphi)id_{V_i}\ \text{for\ all}\ \varphi\in\text{End}_{k}(V_i).\end{equation}
Such an element $c_i$, depending only on the isomorphism class of $V_i$, is called the Schur element associated to $V_i$ (see \cite[Theorem 7.2.1]{GP}). Since $H$ is semisimple, it follows from \cite[Theorem 7.2.6]{GP} that the Schur element $c_i\neq0$ in $\kk$
and the Frobenius homomorphism $\lambda\leftharpoonup u$  can be written explicitly as follows: \begin{equation}\label{equ4}\lambda\leftharpoonup u=u\rightharpoonup\lambda=\sum_{i\in I}\frac{1}{c_i}\chi_i.\end{equation}

\section{\bf A formula for the square of antipodes}
In this section, we will provide a formula for $S^2$ by virtue of a nonzero integral $\Lambda$ of $H$. Then we study some properties of the element  $\mathbf{u}:=S(\Lambda_{(2)})\Lambda_{(1)}$. Especially, we will give a sufficient and necessary condition for $S^2=id$ via the integral $\Lambda$.

Let $u$ be a unit of $H$ satisfying  $S^2(a)=uau^{-1}$ for all $a\in H$. We fix a left integral $\Lambda$ of $H$ and a right integral $\lambda$ of $H^*$ such that $\lambda(\Lambda)=1$. We denote $\{V_i\}_{i\in I}$ the set of all simple left $H$-modules up to isomorphism. For each $V_i$ we denote $c_i$ the Schur element of $V_i$ associated to the dual basis $\Lambda_{(1)}\otimes u^{-1}S(\Lambda_{(2)})$ of $H$ with respect to the Frobenius homomorphism $\lambda\leftharpoonup u$. We denote $\{e_i\}_{i\in I}$ the set of all central primitive idempotents of $H$.
We first establish a relationship between the elements $u$ and $\mathbf{u}=S(\Lambda_{(2)})\Lambda_{(1)}$.

\begin{prop}\label{th1} With the notions above, we have
$\mathbf{u}=u\sum_{i\in I}\dim_{\kk}(V_i)c_ie_i$, which is a unit of $H$.
\end{prop}
\proof Note that each central primitive idempotent $e_i$ acts as the identity on $V_i$ and annihilates $V_j$ for $j\neq i$. It follows that $\chi_j(e_i)=\dim_{\kk}(V_i)$ if $i=j$ and 0 otherwise. By (\ref{equ4}) we have
$$\chi_i(a)=\chi_i(ae_i)=\sum_{j\in I}\frac{1}{c_j}\chi_j(c_iae_i)=(u\rightharpoonup\lambda)(c_iae_i)=(uc_ie_i\rightharpoonup\lambda)(a).$$ Thus, $\chi_i=uc_ie_i\rightharpoonup\lambda$ and hence \begin{equation}\label{equ3}\chi_{H}=\sum_{i\in I}\dim_{\kk}(V_i)\chi_i=u\sum_{i\in I}\dim_{\kk}(V_i)c_ie_i\rightharpoonup\lambda.\end{equation}
For any map $\varphi\in\text{End}_{\kk}(H)$, the trace of $\varphi$ is $\text{tr}(\varphi)=\lambda(\varphi(S(\Lambda_{(2)}))\Lambda_{(1)})$ (see \cite[Theorem 2]{Rad}). Taking into account that $\varphi=L_a$, where $L_a$ is the left multiplication operator of $H$ by $a$, we have
$$\chi_H(a)=\text{tr}(L_a)=\lambda(aS(\Lambda_{(2)})\Lambda_{(1)})=(S(\Lambda_{(2)})\Lambda_{(1)}\rightharpoonup\lambda)(a).$$ This implies that $\chi_H=S(\Lambda_{(2)})\Lambda_{(1)}\rightharpoonup\lambda.$ Comparing it with (\ref{equ3}) and using the non-degeneracy of the Frobenius homomorphism $\lambda$, we have $$S(\Lambda_{(2)})\Lambda_{(1)}=u\sum_{i\in I}\dim_{\kk}(V_i)c_ie_i.$$ Since $p>\dim_{\kk}(H)^{1/2}$, it follows that $p^2>\dim_{\kk}(H)=\sum_{i\in I}\dim_{\kk}(V_i)^2\geq\dim_{\kk}(V_i)^2$. Hence $p>\dim_{\kk}(V_i)$ and $\dim_{\kk}(V_i)\neq0$ in $\kk$ for any $i\in I$.
Thus, the element $u$ is the same as $S(\Lambda_{(2)})\Lambda_{(1)}$ up to a central unit $\sum_{i\in I}\dim_{\kk}(V_i)c_ie_i$.
\qed

\begin{rem} Proposition \ref{th1} also holds if the field $\kk$ has characteristic 0. In this case,  $S^2=id$ (see \cite{LR1} or \cite{LR}) implying that  $\mathbf{u}=S(\Lambda_{(2)})\Lambda_{(1)}=S(\Lambda_{(2)})S^2(\Lambda_{(1)})=S(S(\Lambda_{(1)})\Lambda_{(2)})=\varepsilon(\Lambda)$.\end{rem}

Proposition \ref{th1} gives a formula for $S^2$, namely,  $$S^2(a)=\mathbf{u}a\mathbf{u}^{-1}\ \text{for}\ a\in H,$$ where $\mathbf{u}=S(\Lambda_{(2)})\Lambda_{(1)}.$
In the sequel, we will replace $u$ with $\mathbf{u}$. In this case, the equality (\ref{equ5}) turns out to be \begin{equation}\label{equ55}\Lambda_{(1)}\otimes \mathbf{u}^{-1}S(\Lambda_{(2)})=\mathbf{u}^{-1}S(\Lambda_{(2)})\otimes \Lambda_{(1)},\end{equation}
which is the dual basis of $H$ with respect to the Frobenius homomorphism $\lambda\leftharpoonup \mathbf{u}$.
The Schur element associated to the simple $H$-module $V_i$ under the new dual basis $\Lambda_{(1)}\otimes \mathbf{u}^{-1}S(\Lambda_{(2)})$ with respect to the Frobenius homomorphism $\lambda\leftharpoonup \mathbf{u}$ is $\frac{1}{\dim_{\kk}(V_i)}$.
Therefore, the equality (\ref{equ4}) turns out to be
\begin{equation}\label{equ0077}\lambda\leftharpoonup \mathbf{u}=\mathbf{u}\rightharpoonup\lambda=\sum_{i\in I}\dim_{\kk}(V_i)\chi_i=\chi_H.\end{equation}
By applying \cite[Theorem 1.5]{Doi} and (\ref{equ55}), we obtain the expression of
each central primitive idempotent $e_i$ of $H$ as follows:
\begin{equation}\label{equ8-112}e_i=\dim_{\kk}(V_i)\chi_i(\Lambda_{(1)})\mathbf{u}^{-1}S(\Lambda_{(2)})=\dim_{\kk}(V_i) \chi_i(\mathbf{u}^{-1}S(\Lambda_{(2)}))\Lambda_{(1)}.\end{equation}

Let $g\in G(H)$ and $\alpha\in\text{Alg}(H,k)$ be the modular elements of $H$ and $H^*$ respectively. Recall that the Radford's formula of $S^4$ has the form (see \cite[Proposition 6]{Rad1}): $$S^4(a)=\alpha^{-1}\rightharpoonup(gag^{-1})\leftharpoonup\alpha.$$ Since $H$ is unimodular, i.e., $\alpha=\varepsilon$, the Radford's formula of $S^4$ now becomes $$S^4(a)=gag^{-1}.$$
The distinguished group-like element $g$ and the integral $\Lambda$ of $H$ satisfy the following useful equality (see \cite[Theorem 3(d)]{Rad}):
\begin{equation}\label{equ0}\Lambda_{(2)}\otimes\Lambda_{(1)}=\Lambda_{(1)}\otimes S^2(\Lambda_{(2)})g.\end{equation}
After these preparations, we give some properties of the element $\mathbf{u}$ as follows:

\begin{prop}\label{prop3.20}
The element $\mathbf{u}=S(\Lambda_{(2)})\Lambda_{(1)}$ satisfies the the following properties:
\begin{enumerate}
  \item $\mathbf{u}=\chi_{H}(\Lambda_{(1)})S(\Lambda_{(2)})$.
  \item $\Lambda_{(1)}\mathbf{u}^{-1} S(\Lambda_{(2)})=1.$
  \item $\lambda(e_i)=\dim_{\kk}(V_i)\chi_i(\mathbf{u}^{-1})$.
  \item $\mathbf{u}S(\mathbf{u})=S(\mathbf{u})\mathbf{u}=\varepsilon(\Lambda)\sum_{i\in I}\frac{\dim_{\kk}(V_i)^2}{\lambda(e_i)}e_i$.
  \item $S(\mathbf{u}^{-1})\mathbf{u}=\mathbf{u}S(\mathbf{u}^{-1})$, which is the distinguished group-like element $g$ of $H$.
\end{enumerate}
\end{prop}
\proof (1) It follows from (\ref{equ8-112}) that $e_i\mathbf{u}=\dim_{\kk}(V_i)\chi_i(\Lambda_{(1)})S(\Lambda_{(2)})$. Thus, $$\mathbf{u}=\sum_{i\in I}e_i\mathbf{u}=\sum_{i\in I}\dim_{\kk}(V_i)\chi_i(\Lambda_{(1)})S(\Lambda_{(2)})=\chi_{H}(\Lambda_{(1)})S(\Lambda_{(2)}).$$

(2) Since $\Lambda_{(1)}\otimes \mathbf{u}^{-1}S(\Lambda_{(2)})=\mathbf{u}^{-1}S(\Lambda_{(2)})\otimes \Lambda_{(1)}$ by (\ref{equ55}), we obtain the desired result by multiplying the tensor factors together.

(3) Since $e_i=\dim_{\kk}(V_i)\chi_i(\Lambda_{(1)})\mathbf{u}^{-1}S(\Lambda_{(2)})$, it follows that $$e_i=\mathbf{u}e_i\mathbf{u}^{-1}=\dim_{\kk}(V_i)\chi_i(\Lambda_{(1)})S(\Lambda_{(2)})\mathbf{u}^{-1}.$$
Hence $$\lambda(e_i)=\dim_{\kk}(V_i)\chi_i(\Lambda_{(1)})\lambda(S(\Lambda_{(2)})\mathbf{u}^{-1})=\dim_{\kk}(V_i)\chi_i(\mathbf{u}^{-1}),$$
where the last  equality follows from (\ref{equ01}).

(4) For any $a\in H$, we have $S^3(a)=S(S^2(a))=S(\mathbf{u} a\mathbf{u}^{-1})=S(\mathbf{u}^{-1})S(a)S(\mathbf{u})$, we also have $S^3(a)=S^2(S(a))=\mathbf{u}S(a)\mathbf{u}^{-1}$. It follows that $S(\mathbf{u})\mathbf{u}$ is a central unit of $H$. The equality $\mathbf{u}S(\mathbf{u})=S(\mathbf{u})\mathbf{u}$ holds because $S(\mathbf{u})=S(S^2(\mathbf{u}))=S^2(S(\mathbf{u}))=\mathbf{u}S(\mathbf{u})\mathbf{u}^{-1}$.
For the central unit $\mathbf{u}S(\mathbf{u})$, we suppose that $\mathbf{u}S(\mathbf{u})=\sum_{i\in I} k_ie_i$, where each scalar $k_i\neq0$ in $\kk$. Then $e_i\mathbf{u}^{-1}=\frac{1}{k_i}e_iS(\mathbf{u})$. We have
\begin{align*}\lambda(e_i)&=(\mathbf{u}^{-1}\rightharpoonup\chi_H)(e_i)=\chi_H(e_i\mathbf{u}^{-1})=\frac{1}{k_i}\chi_H(e_iS(\mathbf{u}))\\
&=\frac{\dim_{\kk}(V_i)}{k_i}\chi_i(e_iS(\mathbf{u}))=\frac{\dim_{\kk}(V_i)}{k_i}\chi_i(S(\mathbf{u}))\\
&=\frac{\dim_{\kk}(V_i)}{k_i}(\chi_i\circ S)(\mathbf{u})=\frac{\dim_{\kk}(V_i)}{k_i}(\chi_i\circ S)( S(\Lambda_{(2)})\Lambda_{(1)})\\
&=\frac{\dim_{\kk}(V_i)}{k_i}(\chi_i\circ S)(\Lambda_{(1)}S(\Lambda_{(2)}))=\frac{\dim_{\kk}(V_i)^2\varepsilon(\Lambda)}{k_i}\neq0.
\end{align*}
It follows that $k_i=\frac{\dim_{\kk}(V_i)^2\varepsilon(\Lambda)}{\lambda(e_i)}$ and $\mathbf{u}S(\mathbf{u})=\sum_{i\in I} k_ie_i=\varepsilon(\Lambda)\sum_{i\in I} \frac{\dim_{\kk}(V_i)^2}{\lambda(e_i)}e_i$.

(5) Note that $\Lambda_{(2)}\otimes\Lambda_{(1)}=\Lambda_{(1)}\otimes S^2(\Lambda_{(2)})g$ by (\ref{equ0}). Applying $S\otimes id$ to both sides of this equality and multiplying the tensor factors together, we have $\mathbf{u}=S(\mathbf{u})g$ or $g=S(\mathbf{u}^{-1})\mathbf{u}$.
\qed

As a consequence, we obtain the following result:

\begin{cor}\label{cor4}For any central primitive idempotent $e_i$ of $H$,
we have $\lambda(e_i)=\lambda(S(e_i))$.
\end{cor}
\proof  We denote $S(e_i)=e_{i^*}$ for some $i^*\in I$, then  $V_i^*\cong V_{i^*}$, or equivalently, $\chi_i\circ S=\chi_{i^*}$ (see \cite[Lemma 1.8]{Doi}).
By Proposition \ref{prop3.20} (3) we have $$\lambda(S(e_i))=\lambda(e_{i^*})=\dim_{\kk}(V_{i^*})\chi_{i^*}(\mathbf{u}^{-1})=\dim_{\kk}(V_i)\chi_i(S(\mathbf{u}^{-1})).$$
Since $\mathbf{u}S(\mathbf{u})=\varepsilon(\Lambda)\sum_{i\in I} \frac{\dim_{\kk}(V_i)^2}{\lambda(e_i)}e_i$, it follows that $S(\mathbf{u}^{-1})=\mathbf{u}\frac{1}{\varepsilon(\Lambda)}\sum_{i\in I} \frac{\lambda(e_i)}{\dim_{\kk}(V_i)^2}e_i$. Thus,
\begin{align*}\lambda(S(e_i))&=\dim_{\kk}(V_i)\chi_i(S(\mathbf{u}^{-1}))=\frac{\lambda(e_i)}{\varepsilon(\Lambda)\dim_{\kk}(V_i)}\chi_i(\mathbf{u})\\
&=\frac{\lambda(e_i)}{\varepsilon(\Lambda)\dim_{\kk}(V_i)}\chi_i(\Lambda_{(1)}S(\Lambda_{(2)}))=\lambda(e_i).
\end{align*}
We complete the proof. \qed

If the field $\kk$ has characteristic 0, then the antipode $S$ of $H$ satisfies $S^2=id$ (see \cite{LR1} or \cite{LR}). This further implies that the integral $\Lambda$ of $H$ is cocommutative (see \cite[Proposition 2(b)]{LR1}). The following result shows that $\Lambda$ being cocommutative is equivalent to $S^2=id$ when the characteristic of the field $\kk$ is larger then $\dim_{\kk}(H)^{1/2}$.

\begin{prop}\label{thm4}Let $H$ be a finite dimensional semisimple Hopf algebra over the field $\kk$ of characteristic $p>\dim_{\kk}(H)^{1/2}$. The following statements are equivalent:
\begin{enumerate}
  \item The nonzero integral $\Lambda$ of $H$ is cocommutative.
  \item The nonzero integral $\lambda$ of $H^*$ is cocommutative.
  \item $S^2=id$.
\end{enumerate}
\end{prop}
\proof It can be seen from  \cite[Corollary 5]{Rad} that Part (2) and Part (3) are equivalent. We next show that Part (1) and Part (3) are equivalent. If $\Lambda$ is cocommutative, then $\mathbf{u}=S(\Lambda_{(2)})\Lambda_{(1)}= S(\Lambda_{(1)})\Lambda_{(2)}=\varepsilon(\Lambda)$. It follows from $S^2(a)=\mathbf{u}a\mathbf{u}^{-1}$ that $S^2=id$. Conversely, if $S^2=id$, then $\mathbf{u}=S(\Lambda_{(2)})\Lambda_{(1)}=S(\Lambda_{(2)})S^2(\Lambda_{(1)})=S( S(\Lambda_{(1)})\Lambda_{(2)})=\varepsilon(\Lambda)$. By Proposition \ref{prop3.20}, we have $g=S(\mathbf{u}^{-1})\mathbf{u}=1$. Since $\Lambda_{(2)}\otimes\Lambda_{(1)}=\Lambda_{(1)}\otimes S^2(\Lambda_{(2)})g$ by (\ref{equ0}), it follows that $\Lambda_{(2)}\otimes\Lambda_{(1)}=\Lambda_{(1)}\otimes \Lambda_{(2)}$. We complete the proof.
\qed

\section{\bf Higher FS indicators}
If the field $\kk$ has characteristic 0, the $n$-th FS indicators of finite dimensional representations of semisimple Hopf algebras have been studied in \cite{KSZ}. In this section, we will generalize these indicators from characteristic 0 to the case of characteristic $p>\dim_{\kk}(H)^{1/2}$ and describe them via a nonzero integral $\Lambda$ of $H$. We begin with the following preparations.

Let $H$ be a finite dimensional semisimple Hopf algebra over the field $\kk$ of characteristic $p>\dim_{\kk}(H)^{1/2}$ with a nonzero integral $\Lambda$ and $\mathbf{u}=S(\Lambda_{(2)})\Lambda_{(1)}$.
Applying $\Delta_{n-1}\otimes id$ to both sides of the equality: $\Lambda_{(2)}\otimes\Lambda_{(1)}=\Lambda_{(1)}\otimes S^2(\Lambda_{(2)})g$ (see (\ref{equ0})), we have
\begin{equation*}\Lambda_{(2)}\otimes\cdots\otimes\Lambda_{(n)}\otimes\Lambda_{(1)}=\Lambda_{(1)}\otimes\cdots \otimes\Lambda_{(n-1)} \otimes S^2(\Lambda_{(n)})g.\end{equation*}
Since $g=\mathbf{u}S(\mathbf{u}^{-1})$ and $S^2(\Lambda_{(n)})=\mathbf{u}\Lambda_{(n)}\mathbf{u}^{-1}$, the above equality induces the following equality:
\begin{equation}\label{equ00}\Lambda_{(2)}\otimes\cdots\otimes\Lambda_{(n)}\otimes \mathbf{u}^{-1}\Lambda_{(1)}=\Lambda_{(1)}\otimes\cdots \otimes\Lambda_{(n-1)} \otimes \Lambda_{(n)}S(\mathbf{u}^{-1}).\end{equation}

Note that the category Rep$(H)$ of finite dimensional representations of $H$ is a semisimple tensor category. Let $j_{\mathbf{u}}:id\rightarrow(-)^{**}$ be a natural isomorphism between the identity functor and
the functor of taking the second dual. It is completely determined by a collection of $H$-module isomorphisms
$$j_{\mathbf{u},V}:V\rightarrow V^{**},\ \ j_{\mathbf{u},V}(v)(f)=f(\mathbf{u}v)\ \text{for}\ v\in V, f\in V^*.$$
The inverse of $j_{\mathbf{u},V}$ is $$j_{\mathbf{u},V}^{-1}:V^{**}\rightarrow V,\ \ \alpha\mapsto j_{\mathbf{u},V}^{-1}(\alpha),$$ where $j_{\mathbf{u},V}^{-1}(\alpha)\in V$ satisfies the equality  $f(j_{\mathbf{u},V}^{-1}(\alpha))=\alpha(S^{-1}(\mathbf{u}^{-1})f)$ for $f\in V^*$. Since $S^2(h)=\mathbf{u}h\mathbf{u}^{-1}$ and $\mathbf{u}$ is not known to be a group-like element, the natural isomorphism $j_{\mathbf{u}}$ is not necessarily a tensor isomorphism.
Although the representation category Rep$(H)$ with respect to the structure $j_{\mathbf{u}}$ is not necessarily pivotal, we may still define higher FS indicators for any finite dimensional representation of $H$  using the structure $j_{\mathbf{u}}$ of Rep$(H)$.

We denote $V^{\otimes n}$ the $n$-th tensor power of $V$ where $V^{\otimes0}$ is the trivial $H$-module $\kk$.
For any natural number $n\geq1$,
we define the following $\kk$-linear map $$E_V^{n}:\text{Hom}_H(\kk,V^{\otimes n})\rightarrow \text{Hom}_H(\kk,V^{\otimes n}),\ \ f\mapsto E_V^{n}(f),$$
where $E_V^{n}(f)$ is an $H$-module morphism from $\kk$ to $V^{\otimes n}$ given by
$$E_V^{n}(f):\kk\xrightarrow{\text{coev}_{V^*}}V^*\otimes V^{**}=V^*\otimes \kk\otimes V^{**}\xrightarrow{id\otimes f\otimes id}
V^*\otimes V^{\otimes n}\otimes V^{**}$$
$$\xrightarrow{\text{ev}_V\otimes id}V^{\otimes (n-1)}\otimes V^{**}\xrightarrow{id\otimes j_{\mathbf{u},V}^{-1}}V^{\otimes n}.$$
Here the maps $\text{coev}_{V^*}$ and $\text{ev}_V$ are the usual coevaluation morphism of $V^*$ and evaluation morphism of $V$ respectively.
If we set $f(1)=\sum v_1\otimes\cdots\otimes v_n\in V^{\otimes n}$, the above definition of $E_V^{n}(f)$ shows that
\begin{equation}\label{equ09}E_V^{n}(f)(1)=\sum v_2\otimes\cdots\otimes v_n\otimes \mathbf{u}^{-1}v_1.\end{equation}

Similar to \cite{NS2}, we give the definition of the $n$-th FS indicator of $V$ to be the trace of the linear operator $E_V^{n}$ as follows:

\begin{defi}Let $H$ be a finite dimensional semisimple Hopf algebra over the field $\kk$ of characteristic $p>\dim_{\kk}(H)^{1/2}$.
For any finite dimensional representation $V$ of $H$, the $n$-th FS indicator of $V$ is defined by
$$\nu_n(V)=\text{tr}(E_V^{n})\ \text{for}\ n\geq1.$$
\end{defi}

Similar to the characteristic 0 case, the $n$-th FS indicator of $V$ defined above can also be described by a nonzero integral $\Lambda$ of $H$:

\begin{thm}\label{prop3}Let $\Lambda$ be a nonzero integral of $H$  and $\mathbf{u}=S(\Lambda_{(2)})\Lambda_{(1)}$. Suppose $\chi_V$ is the character of a finite dimensional representation $V$ of $H$.
We have $$\nu_n(V)=\chi_V(\mathbf{u}^{-1}\Lambda_{(1)}\cdots\Lambda_{(n)})\ \text{for}\ n\geq1.$$
\end{thm}
\proof We first show that the equality $\nu_n(V)=\chi_V(\mathbf{u}^{-1}\Lambda_{(1)}\cdots\Lambda_{(n)})$ holds for an idempotent integral $\Lambda$.
Suppose that $\alpha$ is the following $\kk$-linear map
$$\alpha:V^{\otimes n}\rightarrow V^{\otimes n},\ \ v_1\otimes v_2\otimes\cdots \otimes v_n\mapsto v_2\otimes\cdots \otimes v_n\otimes v_1$$ and $\delta=\alpha\circ(\mathbf{u}^{-1}\Lambda_{(1)}\otimes \Lambda_{(2)}\otimes\cdots\otimes\Lambda_{(n)}).$ We have
\begin{align}\label{equ29}\nonumber\delta(v_1\otimes v_2\otimes\cdots \otimes v_n)&=\alpha(\mathbf{u}^{-1}\Lambda_{(1)}v_1\otimes \Lambda_{(2)}v_2\otimes \cdots\otimes\Lambda_{(n)}v_n)\\
\nonumber&=\Lambda_{(2)}v_2\otimes\cdots\otimes \Lambda_{(n)}v_n\otimes \mathbf{u}^{-1}\Lambda_{(1)}v_1\\
&=\Lambda_{(1)}v_2\otimes\cdots\otimes \Lambda_{(n-1)}v_n\otimes \Lambda_{(n)}S(\mathbf{u}^{-1})v_1\ \text{by}\ (\ref{equ00})\\
\nonumber&=\Lambda\cdot(v_2\otimes\cdots\otimes v_n\otimes S(\mathbf{u}^{-1})v_1).
\end{align}
This shows that $\delta(V^{\otimes n})\subseteq \Lambda\cdot V^{\otimes n}=(V^{\otimes n})^H.$
Note that the map $$\Phi:\textrm{Hom}_H(\kk,V^{\otimes n})\rightarrow(V^{\otimes n})^H,\ \ f\mapsto f(1)$$ is an $H$-module isomorphism. We claim that the following diagram is commutative:
\begin{linenomath}
\begin{equation*}
\xymatrix{
  \textrm{Hom}_H(\kk,V^{\otimes n}) \ar[d]_{\Phi} \ar[r]^{E_V^{n}}& \textrm{Hom}_H(\kk,V^{\otimes n})\ar[d]^{\Phi}& \\
 (V^{\otimes n})^H \ar[r]^{\delta}&(V^{\otimes n})^H.&}
\end{equation*}
\end{linenomath}
Indeed, for any $f\in\text{Hom}_H(\kk,V^{\otimes n})$, we suppose that $f(1)=\sum v_1\otimes\cdots\otimes v_n\in V^{\otimes n}$. It follows from $f(1)=f(\Lambda\cdot1)=\Lambda\cdot f(1)$ that
\begin{equation}\label{equ19}\sum v_1\otimes\cdots\otimes v_n=\sum\Lambda_{(1)}v_1\otimes\cdots\otimes\Lambda_{(n)}v_n.\end{equation}
On the one hand, we have
\begin{align*}
(\delta\circ\Phi)(f)&=\delta(f(1))=\delta(\sum v_1\otimes\cdots\otimes v_n)\\
&=\Lambda\cdot(\sum v_2\otimes\cdots\otimes v_n\otimes S(\mathbf{u}^{-1})v_1)\ \ \text{by}\ (\ref{equ29})
\end{align*}
On the other hand, we have
\begin{align*}
(\Phi\circ E_V^n)(f)&=E_V^n(f)(1)=\sum v_2\otimes\cdots\otimes v_n\otimes \mathbf{u}^{-1}v_1\ \ \ \text{by}\ (\ref{equ09})\\
&=\sum \Lambda_{(2)}v_2\otimes\cdots\otimes \Lambda_{(n)}v_n\otimes \mathbf{u}^{-1}\Lambda_{(1)}v_1\ \ \ \text{by}\ (\ref{equ19})\\
&=\sum \Lambda_{(1)}v_2\otimes\cdots\otimes \Lambda_{(n-1)}v_n\otimes\Lambda_{(n)}S(\mathbf{u}^{-1})v_1\ \ \ \text{by}\ (\ref{equ00})\\
&=\Lambda\cdot(\sum v_2\otimes\cdots\otimes v_n\otimes S(\mathbf{u}^{-1})v_1).
\end{align*}
We obtain that $\delta\circ\Phi=\Phi\circ E_V^n$, or equivalently, $E_V^n=\Phi^{-1}\circ\delta\circ\Phi$. It follows that
\begin{align*}\nu_n(V)&=\text{tr}(E_V^{n})=\text{tr}_{V^{\otimes n}}(\delta)\\
&=\text{tr}_{V^{\otimes n}}(\alpha\circ(\mathbf{u}^{-1}\Lambda_{(1)}\otimes \Lambda_{(2)}\otimes\cdots\otimes\Lambda_{(n)}))\\
&=\text{tr}_V(\mathbf{u}^{-1}\Lambda_{(1)}\cdots\Lambda_{(n)})\\
&=\chi_V(\mathbf{u}^{-1}\Lambda_{(1)}\cdots\Lambda_{(n)}),
\end{align*}
where the equality $\text{tr}_{V^{\otimes n}}(\alpha\circ(\mathbf{u}^{-1}\Lambda_{(1)}\otimes \Lambda_{(2)}\otimes\cdots\otimes\Lambda_{(n)}))
=\text{tr}_V(\mathbf{u}^{-1}\Lambda_{(1)}\cdots\Lambda_{(n)})$ follows from \cite[Lemma 2.3]{KSZ}. We have shown that  $\nu_n(V)=\chi_V(\mathbf{u}^{-1}\Lambda_{(1)}\cdots\Lambda_{(n)})$ where $\Lambda$ is idempotent. Since $\mathbf{u}^{-1}\Lambda_{(1)}\cdots\Lambda_{(n)}$ does not depend on the choice of the nonzero integral $\Lambda$, the equality $\nu_n(V)=\chi_V(\mathbf{u}^{-1}\Lambda_{(1)}\cdots\Lambda_{(n)})$ holds for any nonzero integral $\Lambda$ of $H$.
\qed

\begin{rem}
If the field $\kk$ has characteristic 0 and $\Lambda$ is idempotent, then $\mathbf{u}=\varepsilon(\Lambda)=1$. In this case, the $n$-th FS indicator of $V$ is $\chi_V(\Lambda_{(1)}\cdots\Lambda_{(n)})$, which is  the one defined in \cite[Definition 2.3]{KSZ}.
\end{rem}

In the rest of this section, we will extend the $n$-th FS indicator $\nu_n(V)$ of $V$ from $n\geq1$ to the case $n\in\mathbb{Z}$.
Recall that the $n$-th Sweedler power map $P_n:H\rightarrow H$ is defined by
$$P_n(a)=\left\{
           \begin{array}{ll}
             a_{(1)}\cdots a_{(n)}, &  n\geq1; \\
             \varepsilon(a), & n=0; \\
             S(a_{(1)})\cdots S(a_{(-n)}), &  n\leq-1.
           \end{array}
         \right.
$$ From the $n$-th Sweedler power map $P_n$ of $H$, we may see that $$\nu_n(V)=\chi_V(\mathbf{u}^{-1}P_n(\Lambda))\ \text{for}\ n\geq1.$$
However, this expression is well-defined for any integer $n$. Thus, we may extend this formula from $n\geq1$ to any integer $n$ stated as follows:
\begin{defi}Let $H$ be a finite dimensional semisimple Hopf algebra over the field $\kk$ of characteristic $p>\dim_{\kk}(H)^{1/2}$.
For any finite dimensional representation $V$ of $H$ and any $n\in\mathbb{Z}$, the $n$-th FS indicator of $V$ is defined by
$$\nu_n(V)=\chi_V(\mathbf{u}^{-1}P_n(\Lambda)),$$
where $\mathbf{u}=S(\Lambda_{(2)})\Lambda_{(1)}$.
\end{defi}

\begin{rem}\label{rem1}
\begin{enumerate}
\item Note that $S(\Lambda)=\Lambda$. The $n$-th FS indicator of $V$  can be written as
$$\nu_n(V)=\left\{
             \begin{array}{ll}
               \chi_V(\mathbf{u}^{-1}\Lambda_{(1)}\cdots\Lambda_{(n)}), & n\geq1; \\
               \chi_V(\mathbf{u}^{-1}\varepsilon(\Lambda)), & n=0; \\
               \chi_V(\mathbf{u}^{-1}\Lambda_{(-n)}\cdots\Lambda_{(1)}), & n\leq-1.
             \end{array}
           \right.
$$
\item By Proposition \ref{prop3.20} (4), we have \begin{equation*}\mathbf{u}^{-1}S(\mathbf{u}^{-1})=\sum_{i\in I}\frac{\lambda(e_i)}{\varepsilon(\Lambda)\dim_{\kk}(V_i)^2}e_i\in Z(H).\end{equation*} It follows that  \begin{align*}\nu_0(V)&=\varepsilon(\Lambda)\chi_V(\mathbf{u}^{-1})=\varepsilon(\Lambda)\chi_V(\mathbf{u}^{-1}S(\mathbf{u}^{-1})S(\mathbf{u}))\\
&=\sum_{i\in I}\frac{\lambda(e_i)}{\dim_{\kk}(V_i)^2}\chi_V(e_iS(\mathbf{u}))=\sum_{i\in I}\frac{\lambda(e_i)}{\dim_{\kk}(V_i)^2}\chi_V(e_iS(\Lambda_{(1)})S^2(\Lambda_{(2)}))\\
&=\sum_{i\in I}\frac{\lambda(e_i)}{\dim_{\kk}(V_i)^2}\chi_V(e_iS^2(\Lambda_{(2)})S(\Lambda_{(1)}))=\varepsilon(\Lambda)\sum_{i\in I}\frac{\lambda(e_i)}{\dim_{\kk}(V_i)^2}\chi_V(e_i).
\end{align*}
  \item $\nu_{-1}(V)=\nu_{1}(V)=\chi_V(\mathbf{u}^{-1}\Lambda)=\chi_V(\frac{\Lambda}{\varepsilon(\Lambda)})$.
  \item By \cite[Proposition 3.1]{WLL}, $\Lambda_{(1)}\Lambda_{(2)}$ and $\Lambda_{(2)}\Lambda_{(1)}$ are both central elements of $H$, they are determined by the values that the characters $\chi_i$ for all $i\in I$ take on them. It follows from $\chi_i(\Lambda_{(1)}\Lambda_{(2)})=\chi_i(\Lambda_{(2)}\Lambda_{(1)})$ that $\Lambda_{(1)}\Lambda_{(2)}=\Lambda_{(2)}\Lambda_{(1)}$. Therefore, $\nu_{-2}(V)=\nu_{2}(V)$.
\end{enumerate}
\end{rem}

The higher FS indicators of any simple module $V_i$ can be described as follows:
\begin{prop}For any $n\in\mathbb{Z}$ and any simple module $V_i$ with the character $\chi_i$, we have
$$\nu_n(V_i)=\frac{\chi_i(P_n(\Lambda))\lambda(e_i)}{\dim_{\kk}(V_i)^2}.
$$
\end{prop}
\proof
Since $P_n(\Lambda)\in Z(H)$ for any $n\in\mathbb{Z}$ (see \cite[Proposition 3.1]{WLL}), it follows that $P_n(\Lambda)=\sum_{i\in I}\frac{\chi_i(P_n(\Lambda))}{\dim_{\kk}(V_i)}e_i$. The $n$-th FS indicator of  $V_i$ is $$\nu_n(V_i)=\chi_i(\mathbf{u}^{-1}P_n(\Lambda))
=\frac{\chi_i(P_n(\Lambda))}{\dim_{\kk}(V_i)}\chi_i(\mathbf{u}^{-1})
=\frac{\chi_i(P_n(\Lambda))\lambda(e_i)}{\dim_{\kk}(V_i)^2},
$$
where the last equality follows from Proposition \ref{prop3.20} (3).
\qed

For any semisimple Hopf algebra over a field $\kk$ of characteristic 0, the finite dimensional representation $V$ and its dual $V^*$ have the same $n$-th FS indicators for all $n\geq1$ (see \cite[Section 2.3]{KSZ}). The following result shows that this result also holds for the $n$-th FS indicators defined for the Hopf algebra $H$ over the field $\kk$ of characteristic $p>\dim_{\kk}(H)^{1/2}$.

\begin{prop}\label{prop3.21} Let $H$ be a finite dimensional semisimple Hopf algebra over the field $\kk$ of characteristic $p>\dim_{\kk}(H)^{1/2}$. Let $V$ be a finite dimensional representation of $H$ with the dual $V^*$. We have $\nu_n(V)=\nu_n(V^*)$ for all $n\in\mathbb{Z}$.
\end{prop}
\proof
Since $S(\Lambda)=\Lambda$, we have $S(P_n(\Lambda))=P_n(\Lambda)$ for any $n\in\mathbb{Z}$.
For the case $n\geq1$, the $n$-th FS indicator of $V^*$ is
\begin{align*}\nu_n(V^*)&=(\chi_{V^*})(\mathbf{u}^{-1}P_n(\Lambda))=(\chi_V\circ S)(\mathbf{u}^{-1}P_n(\Lambda))\\
&=\chi_V(\Lambda_{(1)}\cdots \Lambda_{(n)}S(\mathbf{u}^{-1}))=\chi_V(\Lambda_{(2)}\cdots\Lambda_{(n)}\mathbf{u}^{-1}\Lambda_{(1)})\ \ \text{by}\ (\ref{equ00})\\
&=\chi_V(\mathbf{u}^{-1}\Lambda_{(1)}\Lambda_{(2)}\cdots\Lambda_{(n)})=\nu_n(V).
\end{align*}

For the case $n\leq-1$,
the $n$-th FS indicator of $V^*$ is
\begin{align*}\nu_{n}(V^*)&=(\chi_{V^*})(\mathbf{u}^{-1}P_n(\Lambda))=(\chi_V\circ S)(\mathbf{u}^{-1}P_n(\Lambda))\\
&=\chi_V(\Lambda_{(-n)}\cdots \Lambda_{(1)}S(\mathbf{u}^{-1}))=\chi_V(S(\mathbf{u}^{-1})\Lambda_{(-n)}\cdots \Lambda_{(1)})\\
&=\chi_V(S(\mathbf{u}^{-1})\mathbf{u}^{-1}\Lambda_{(1)}S(\mathbf{u})\Lambda_{(-n)}\cdots\Lambda_{(2)})\ \ \text{by}\ (\ref{equ00})\\
&=\chi_V(\Lambda_{(1)}\mathbf{u}^{-1}\Lambda_{(-n)}\cdots\Lambda_{(2)})=\chi_V(\mathbf{u}^{-1}\Lambda_{(-n)}\cdots\Lambda_{(1)})\\
&=\nu_{n}(V).
\end{align*}

For the case $n=0$,
we denote $S(e_i)=e_{i^*}$ for any $i\in I$, then $*$ is a permutation of $I$,  $V_{i^*}\cong V_i^*$ and $\lambda(e_{i^*})=\lambda(e_i)$ by Corollary \ref{cor4}. We have
\begin{align*}\nu_0(V^*)&=\varepsilon(\Lambda)\sum_{i\in I}\frac{\lambda(e_i)}{\dim_{\kk}(V_i)^2}\chi_V(S(e_i))\ \ \text{by\ Remark}\ \ref{rem1} (2)\\
&=\varepsilon(\Lambda)\sum_{i\in I}\frac{\lambda(e_{i^*})}{\dim_{\kk}(V_{i^*})^2}\chi_V(e_{i^*})=\varepsilon(\Lambda)\sum_{i\in I}\frac{\lambda(e_{i})}{\dim_{\kk}(V_{i})^2}\chi_V(e_{i})\\
&=\nu_0(V).
\end{align*}
We complete the proof.
\qed

Kashina-Sommerh$\ddot{a}$user-Zhu has shown in \cite[Proposition 2.5]{KSZ} that the $n$-th FS indictor of the regular representation of a semisimple
Hopf algebra over a field of characteristic 0 can be described as $\text{tr}(S\circ P_{n-1})$ for $n\geq1.$
The following result shows that this formula also holds for the $n$-th FS indicators defined for the Hopf algebra $H$ over the field $\kk$ of characteristic $p>\dim_{\kk}(H)^{1/2}$.

\begin{prop}Let $H$ be a finite dimensional semisimple Hopf algebra over the field $\kk$ of characteristic $p>\dim_{\kk}(H)^{1/2}$. For any $n\in\mathbb{Z}$, the $n$-th FS indictor of the regular representation of $H$ can be written as
 $\nu_n(H)=\text{tr}(S\circ P_{n-1}),$ where $P_{n-1}$ is the $(n-1)$-th Sweedler power map of $H$.
\end{prop}
\proof We choose a left integral $\Lambda$ of $H$ and a right integral $\lambda$ of $H^*$ such that $\lambda(\Lambda)=1$. For any $n\in\mathbb{Z}$, by Radford's trace formula \cite[Theorem 2]{Rad}, we have
\begin{align*}
\text{tr}(S\circ P_{n-1})&=\text{tr}(P_{n-1}\circ S)=\lambda(S(\Lambda_{(2)})(P_{n-1}\circ S)(\Lambda_{(1)}))\\
&=\lambda(S(\Lambda_{(2)})P_{n-1}(S(\Lambda_{(1)})))=\lambda(\Lambda_{(1)}P_{n-1}(\Lambda_{(2)}))\\
&=\lambda(P_n(\Lambda))=\chi_H(\mathbf{u}^{-1}P_n(\Lambda))\ \ \text{by}\ (\ref{equ0077})\\
&=\nu_{n}(H).
\end{align*}
We complete the proof. \qed


\section{\bf Gauge invariants}
In this section, we will show that the $n$-th FS indicator $\nu_n(V)$ defined in Section 4 is a gauge invariant of the tensor category Rep$(H)$ for any $n\in\mathbb{Z}$ and any finite dimensional representation $V$ of the semisimple Hopf algebra $H$ .

Recall from \cite{AEGN} that a  (normalized) twist for semisimple Hopf algebra $H$ is an invertible element $J\in H\otimes H$ that satisfies $(\varepsilon\otimes id)(J)=(id\otimes\varepsilon)(J)=1$ and
$$(\Delta\otimes id)(J)(J\otimes1)=(id\otimes\Delta)(J)(1\otimes J).$$
We write $J=J^{(1)}\otimes J^{(2)}$ and $J^{-1}=J^{-(1)}\otimes J^{-(2)}$, where the summation is understood. 

Given a twist $J$ for $H$ one can define a new Hopf algebra $H^{J}$ with the same algebra structure and counit as $H$, for which the comultiplication $\Delta^J$ and antipode $S^J$ are given respectively by
$$\Delta^J(a)=J^{-1}\Delta(a)J,$$
$$S^J(a)=Q^{-1}_JS(a)Q_J,\ \ \text{for}\ a\in H,$$
where $Q_J=S(J^{(1)})J^{(2)}$, which is an invertible element of $H$ with the inverse $Q_J^{-1}=J^{-(1)}S(J^{-(2)})$.
With the notions above, we have the following result:

\begin{prop}\label{thm5}Let $H$ be a finite dimensional semisimple Hopf algebra over the field $\kk$ of characteristic $p>\dim_{\kk}(H)^{1/2}$ and $V$ a finite dimensional representation of $H$.
The $n$-th FS indicator $\nu_n(V)$ of $V$ is invariant under twisting for any $n\in\mathbb{Z}$.
\end{prop}
\proof
Let $\Lambda$ be a nonzero integral of $H$ and $J$ a normalized twist for $H$. It follows from \cite[Theorem 3.4]{WLL} that $P^J_n(\Lambda)=P_n(\Lambda)$, where $P^J_{n}$ and $P_{n}$ are the $n$-th  Sweedler power maps of $H^J$ and $H$ respectively. Moreover, $P_n(\Lambda)$ is a central element of $H$ (see \cite[Proposition 3.1]{WLL}). Since $\Delta^J(\Lambda)=Q_J^{-1}\Lambda_{(1)}\otimes\Lambda_{(2)}Q_J,$ it follows that
\begin{equation}\label{equ006}\mathbf{u}^J:=S^J(\Lambda_{(2)}Q_J)Q_J^{-1}\Lambda_{(1)}=Q_J^{-1}S(Q_J)\mathbf{u},\end{equation}
where $\mathbf{u}=S(\Lambda_{(2)})\Lambda_{(1)}.$ For $H$-module $V$ with the character $\chi_V$, we denote
$V^J$ the same as $V$ as  $\kk$-linear space but thought of as an $H^J$-module. Then the character of $V^J$ is also $\chi_V$. For any $n\in\mathbb{Z}$, we have
\begin{align*}
\nu_n(V^J)&=\chi_V((\mathbf{u}^J)^{-1}P_n^J(\Lambda))\\&=\chi_V(\mathbf{u}^{-1}S(Q_J^{-1})Q_JP_n^J(\Lambda))\ \ \text{by}\ (\ref{equ006})\\
&=\chi_V(\mathbf{u}^{-1}S(Q_J^{-1})Q_JP_n(\Lambda))\\&=\chi_V(\mathbf{u}^{-1}S^2(J^{-(2)})S(J^{-(1)})S(J^{(1)})J^{(2)}P_n(\Lambda))\\
&=\chi_V(J^{-(2)}\mathbf{u}^{-1}S(J^{-(1)})S(J^{(1)})J^{(2)}P_n(\Lambda))\\&=\chi_V(\mathbf{u}^{-1}S(J^{-(1)})S(J^{(1)})J^{(2)}P_n(\Lambda)J^{-(2)})\\
&=\chi_V(\mathbf{u}^{-1}S(J^{-(1)})S(J^{(1)})J^{(2)}J^{-(2)}P_n(\Lambda))\\&=\chi_V(\mathbf{u}^{-1}P_n(\Lambda))\\
&=\nu_n(V).
\end{align*}
We complete the proof.
\qed

We are now ready to state the main result which says that higher FS indicators are gauge invariants of the tensor category Rep$(H)$.
\begin{thm}
Let $H$ and $H'$ be two finite dimensional semisimple Hopf algebras over the field $\kk$ of characteristic $p>\dim_{\kk}(H)^{1/2}$. If $\mathcal{F}:\text{Rep}(H)\rightarrow\text{Rep}(H')$ is an equivalence of tensor categories, then $\nu_n(V)=\nu_n(\mathcal{F}(V))$ for any $n\in\mathbb{Z}$ and any finite dimensional representation $V$ of $H$.
\end{thm}
\proof
Since the  $\kk$-linear equivalence $\mathcal{F}:\text{Rep}(H)\rightarrow\text{Rep}(H')$
is a tensor equivalence, it follows from \cite[Theorem 2.2]{NS1} that $H$ and $H'$ are gauge equivalent in the sense that there exist
a twist $J$ of $H$  such that $H'$ is isomorphic to $H^J$ as bialgebras. Let $\sigma:H'\rightarrow H^J$ be such an isomorphism. Then $\sigma$ is
automatically a Hopf algebra isomorphism. The isomorphism $\sigma$ induces a $\kk$-linear equivalence $(-)^\sigma:\text{Rep}(H)\rightarrow\text{Rep}(H')$ as follows:  for any finite dimensional $H$-module $V$, $V^{\sigma}=V$ as $\kk$-linear space with the $H'$-module action given by $h'v=\sigma(h')v$ for $h'\in H'$, $v\in V$, and $f^\sigma=f$ for any morphism $f$ in Rep$(H)$.
Moreover, the equivalence  $\mathcal{F}$ is naturally isomorphic to the $\kk$-linear equivalence $(-)^\sigma$ (see \cite[Theorem 1.1]{KMN}).
Therefore, $$\nu_n(\mathcal{F}(V))=\nu_n(V^{\sigma}).$$
Let $\Lambda'$ be a nonzero integral of $H'$ and $S'$ the antipode of $H'$. Note that the map $\sigma: H'\rightarrow H^J$ is a Hopf algebra isomorphism. It follows that $\sigma(\Lambda')=\Lambda$, which is a nonzero integral of $H^J$ and $\sigma(P'_n(\Lambda'))=P^J_n(\Lambda)$, where $P'_n$ and $P^J_n$ are the $n$-th Sweedler power maps of $H'$
and $H^J$ respectively. In particular,
$$\sigma((\mathbf{u}')^{-1}P'_n(\Lambda'))
=(\mathbf{u}^J)^{-1}P^J_n(\Lambda),
$$
where $\mathbf{u}'=S'(\Lambda'_{(2)})\Lambda'_{(1)}$ and $\mathbf{u}^J=S^J(\Lambda_{\langle2\rangle})\Lambda_{\langle1\rangle}$.
We have
\begin{align*}
\nu_n(V^{\sigma})&=\chi_{V^{\sigma}}((\mathbf{u}')^{-1}P'_n(\Lambda'))\\
&=\chi_{V^{J}}(\sigma((\mathbf{u}')^{-1}P'_n(\Lambda')))\\
&=\chi_{V^{J}}((\mathbf{u}^J)^{-1}P^J_n(\Lambda))\\
&=\nu_n(V^J)\\
&=\nu_n(V),
\end{align*}where the last equality follows from Proposition \ref{thm5}.
We conclude that $\nu_n(\mathcal{F}(V))=\nu_n(V)$ for any $n\in\mathbb{Z}$ and any finite dimensional representation $V$ of $H$.  \qed

\section*{Acknowledgement}
The authors would like to thank Professor Juan Cuadra for pointing out a mistake in a previous version of the manuscript.
The first author was supported by Qing Lan Project.
The second author was supported by National Natural Science Foundation of China (Grant No. 11722016). The third author was supported by National Natural Science Foundation of China (Grant No. 11871063).

\end{document}